\documentclass[10pt]{article}
\usepackage{amsfonts,amssymb,amsmath,amscd,amsthm}
\usepackage[all]{xy}

\providecommand{\cprod}[2]{\ensuremath{ #1\rtimes_{#2}\mathbb{Z}}}

\providecommand{\ts}[1]{{\textstyle{#1}}}

\newcommand{\Z}{{\mathbb Z}}
\newcommand{\R}{\mathbb R}
\newcommand{\Q}{{\mathbb Q}}
\newcommand{\T}{{\mathbb T}}

\newcommand{\id}{\mathrm{Id}}

\renewcommand{\d}{\mathrm{d}}

\newcommand{{\frakM}}{\mathfrak{M}}
\newcommand{{\frakS}}{\mathfrak{S}}

\renewcommand{\phi}{\varphi}
\newcommand{\E}{{\Lambda}}

\newtheorem{theorem}{Theorem}[section]
\newtheorem{lemma}[theorem]{Lemma}
\newtheorem{corollary}[theorem]{Corollary}
\newtheorem{proposition}[theorem]{Proposition}
\newtheorem{definition}[theorem]{Definition}
\newtheorem{example}[theorem]{Example}

\begin{document}

\title{Eigenvalues, $K$-theory and  Minimal  Flows}
\author{Benjam\'{\i}n Itz\'a-Ortiz} 
\date{
         }
\maketitle

\begin{abstract}
Let $(Y,T)$ be a minimal suspension flow built over a dynamical system $(X,S)$ and 
with (strictly positive, continuous) ceiling function $f\colon X\to\R$. We show that the  eigenvalues of $(Y,T)$ are contained in the range of a trace on the $K_0$-group of $(X,S)$. Moreover, a trace gives an order isomorphism of a subgroup of  $K_0\left(\cprod{C(X)}{S}\right)$ with the group of eigenvalues of $(Y,S)$. Using this result, we relate the values of $t$ for which the time-$t$ map on  minimal suspension flow is minimal, with the $K$-theory of the base of this suspension.  
\end{abstract}

\section*{Introduction}
It is known that there is a relation between the eigenvalues of a minimal con\-ti\-nuous flow and the values of $t$ for which the time-$t$ map is minimal, cf.\ \cite[4.24.1]{G}. On the other hand, it is known that the group of eigenvalues of a standard suspension built over a minimal group rotation is equal to  the range of the unique normalized trace on the $K_0$-group of this group rotation, see 
e.g.\ \cite{Riedel,Packer,E}. The purpose of this paper is to use the eigenvalues of a minimal suspension flow (with ceiling function not necessary constant) as a link to get a correspondence of the values of $t$ for which the time-$t$ map  is minimal with the $K$-theory of the base of this suspension.

In \cite{Sch}, S.\ Schwartzman made one of the first applications of algebraic topology to topological dynamics by associating to each continuous flow $(Y,T)$ with $T$-invariant Borel probability measure $\nu$, an ``average asymptotic cycle" $A\sb\nu\colon H^1(Y,\Z)\to\R$. It is proved in a Theorem in Section~5 of \cite{Sch} that the range of $A\sb\nu$ contains the eigenvalues of $(Y,T)$. Schwartzman also noted a correspondence between cross sections to $(Y,T)$ and certain (multiplicative) semigroup $C$ in $H^1(Y,T)$ which we will call, following the terminology in \cite{Packer}, the {\it positive Schwartzman cone} of $H^1(Y,\Z)$.
Building upon the work of Schwartzman \cite{Sch}, Connes \cite{C}, and Rieffel \cite{RieffelSME}, J.\ Packer \cite{Packer2} showed that when $(Y,T)$ is the standard suspension flow built over a dynamical system $(X,S)$, the restriction of the Connes' isomorphism $\Phi\colon K_1(C(Y))\to K_0\left(\cprod{C(X)}{S}\right)$ to the positive Schwartzman cone $C$ of $H^1(Y,\Z)$ gives an injective semigroup homomorphism into the positive cone  $K_0\left(\cprod{C(X)}{S}\right)^+$ of $K_0\left(\cprod{C(X)}{S}\right)$.  Furthermore, Packer showed that if $\tau$ is a normalized trace on $\cprod{C(X)}{S}$ then the images of $\tau\circ\Phi$ and $A\sb\nu$ coincide on $C\subset K_1(C(Y))$, and so $\left(\tau\circ\Phi\right)(C)$  contains all strictly positive eigenvalues of $(Y,T)$. This fact made us suspect (something we prove in this work) the existence of    
a  suitable (multiplicative) subgroup $G$  of $H^1(Y,Z)$ such that a trace on $\cprod{C(X)}{S}$ gives an isomorphism of $\Phi( G )\subset K_0\left(\cprod{C(X)}{S}\right)$ with the group of eigenvalues of $(Y,T)$.

We divide this paper in three sections. In Section~1 we introduce notation, terminology and some  results needed in the rest of the paper. We state and prove the main results in Section~2. In the last section we give some examples illustrating our results. 
The idea of a connection between the values of $t$ for which the time-$t$ map on a minimal suspension flow is minimal, and  the $K$-theory of the base of this suspension, was proposed in the author's Ph.D.\ thesis, under the supervision of Professor N.C.\ Phillips.  The author 
gratefully acknowledge support and encouragement from D.\ Handelman, 
T.\ Giordano, and V.\ Pestov. In particular, I am indebted to T.\ Giordano for many stimulating conversations. 


\setcounter{section}{0}
\section{Preliminaries}
In this section we will provide the notation, definitions and results needed in the rest of the paper. Most of the results here are known; in fact, the reader will be referred to their place in the literature whenever we were able to find a reference.

Consider a  compact metric space $Y$ and a continuous map $T\colon Y\times\R\to Y$. Given $(y,r)$ in $Y\times\R$, we choose to write $T(y,t)$ as $T^t(y)$. If  $T^s\circ T^t=T^{s+t}$ for all $s,t\in\R$ and $T^0=\id_Y$ then we say that $(Y,T)$ is a {\em continuous flow}.

A continuous function $\chi\colon Y\to S^1$ is said to be an eigenvector with respect to the flow $(Y,T)$ provided there is a real number $\lambda$ such that
\[
      \chi\left( T^t(y)\right) = e^{2\pi i \lambda t} \chi(y),   \,\,\,\forall y\in Y,\,\,\forall t\in\R.
\]
We call $\lambda$ the eigenvalue associated with $\chi$ or say that $\chi$ is an eigenvector for $\lambda$.

\begin{definition}\label{defOFeigenvalue}
The set of eigenvalues of $(Y,T)$ is denoted by $\E(Y,T)$.
We also set $\E(Y,T)^+=\{\lambda\geq 0 \colon \lambda\in\E(Y,T)  \}$. 
\end{definition}

Let $(Y,T)$ be a continuos flow and let $M$ be a subset of $Y$. We say that $M$ is $T$-invariant if $T^t(M)=M$ for all $t\in\R$, and  we say that $M$ is {\it minimal} (with respect  to $T$) if it is nonempty, closed and $T$-invariant while no proper subset of $M$ has these properties, i.e., $\emptyset$ and $M$ are the only closed $T$-invariant subsets of $M$. When $Y$ itself is a minimal subset of $Y$, then $(Y,T)$ is called a {\it minimal flow}.  Recall that the orbit of $y\in Y$ under $T$ is the set $\{T^t(y)\colon t\in\R\}$. A continuous flow is minimal if and only if every orbit is dense, cf.\ \cite[Lemma~II.3.7]{V}.

\begin{proposition}\label{eigen}
Let $(Y,T)$ be a minimal continuous flow.  
  \begin{enumerate}
      \item If $\chi_1$ and $\chi_2$ are two eigenvectors for the same eigenvalue 
              $\lambda$ of $(Y,T)$ then there exists $c$ in $S^1$ such that $\chi_1=c\chi_2$.
      \item The set $\E(Y,T)$ is a countable subgroup of $\R$;
  \end{enumerate}   
\end{proposition}
\begin{proof}
This is a standard result. A proof can be obtained by imitating the techniques in the proof of \cite[Theorem~5.17]{W}. See also \cite[2.1.5]{G}. Here is a sketch of the proof: For part (1), just observe that $\textstyle\frac{\chi_1}{\chi_2}$ is constant on orbits of points. For part (2), it is clear that $\E(Y,T)$ is a subgroup of $\R$. Let $\chi_1$ and $\chi_2$ be two eigenvectors corresponding to distinct eigenvalues $\lambda_1$ and $\lambda_2$.  We show that $\| \chi_1-\chi_2\| >\ts{\frac{1}{4}}$. Since $C(Y)$ is separable, this will complete the proof. As $\lambda=\lambda_1-\lambda_2$ is a nonzero eigenvalue with eigenvector $\chi_1\overline{\chi_2}$, we may choose $t_0$ in  $\R$ and $y_0$ in $Y$ such that $e^{2\pi i \lambda t_0}\chi_1\overline{\chi_2}(y_0)$ is in the left-hand half of the unit circle. Hence
\begin{eqnarray*}
     \|\chi_1-\chi_2\| & = & \left\|\chi_1\overline{\chi_2} - 1\right\| \\
               & =& \sup_{y\in Y} \left\| \chi_1\overline{\chi_2}\left(T^t(y)\right) - 1 \right\|\\
               & \geq & \left\|e^{2\pi i \lambda t_0}\chi_1\overline{\chi_2}(y_0) -1\right\| \\
               & > & \frac{1}{4}.
\end{eqnarray*} 
\end{proof}

The following observation can be found in \cite[II.3.6(3) and II.10.3(6)]{V}.

\begin{lemma}\label{MinImpliesConnected}
If $(Y,T)$ is a minimal flow then $Y$ is connected. 
\end{lemma}
\begin{proof}
Let $y\in Y$. As the orbit of $y$ is connected then so is the closure of such orbit which, by minimality, is all of $Y$. Thus $Y$ is connected. 
\end{proof}

By a {\it discrete flow} or {\it dynamical system} we mean a pair $(X,S)$ where $X$ is a compact metric space and $S\colon X\to X$ is a homeomorphism.  A discrete flow $(X,S)$ is minimal if there is no nontrivial closed $S$-invariant subset of $X$. Given a continuous flow $(Y,T)$ it induces, for each $t\in\R$, a discrete flow $(Y,T^t)$. We sometimes refer to the map $T^t$ as the time-$t$ map on $Y$. We  now prove the following.

\begin{lemma}\label{min}
Let $(Y,T)$ be a continuous flow where $Y$ is a connected space. Let $t$ be a nonzero real number. The map $T^t$ is not minimal if and only if, for every rational number $r$, the map $T^{rt}$ is not minimal.
\end{lemma}
\begin{proof}
  Suppose that the map $T^t$ is not minimal and let $r$ be a rational number. We show that the map $T^{rt}$ is not minimal.  Assume that $r$ is nonzero. (If $r$ is zero then $T^{rt}=T^0$ is the identity function and so $T^{rt}$  is not minimal.) Write $r={\ts{\frac{p}{q}}}$ with $p$ and $q>0$ relatively prime. We claim that $T^{{\frac{t}{q}}}$ is not minimal.  This will then imply that $T^{rt}=(T^{\frac{t}{q}})^p$ is not minimal (because if $M$ is a nontrivial, closed and $T^{\frac{t}{q}}$-invariant subset of $Y$, then $M$ is also $T^{\frac{pt}{q}}$-invariant), which is what is wanted.
  To prove the claim, assume that  $T^{\frac{t}{q}}$ is minimal. Then, since by hypothesis $T^t=(T^{\frac{t}{q}})^q$ is not minimal,   there is a nontrivial, closed and $T^t$-invariant subset $M$ of $Y$, cf.\ \cite[Theorem~5.2]{W}. Following  \cite[II.9.6(7)]{V}, we see that $M$ is also open, a contradiction, since $Y$ is assumed to be connected.  This proves the claim. The converse of the lemma is obvious.
\end{proof}

If a continous flow $(Y,T)$ is minimal, it is natural to ask for which $t$ the induced discrete flow $(Y,T^t)$ is minimal. The answer is contained in 
the following result which  is essentially outlined in \cite[4.24.1]{G}.

\begin{proposition}\label{step1}
   Let $(Y,T)$ be a minimal continuous flow.  Consider the map with domain $\Q\otimes\E(Y,T)$ and codomain $\R$ defined by $r\otimes\lambda\mapsto r\lambda$ for all $(r,\lambda)$ in $\Q\times\E(Y,T)$. This map is an $\Q$-linear monomorphism with range equal to
    \[
    \Bigl\{0\Bigr\}\cup\left\{t\in\R\setminus\{0\}\colon T^{{\frac{1}{t}}}\text{ is not minimal}\right\}.
   \]
Hence the set above is a countable $\Q$-linear subspace of $\R$ isomorphic to $\Q\otimes \E(Y,T)$.
\end{proposition}
\begin{proof}
  As $\E(Y,T)$ is an abelian (countable) subgroup of $\R$, cf.\ Lemma~\ref{eigen}, one readily  sees that every element in $\Q\otimes\E(Y,T)$ can be written as a pure tensor and that the map  with domain $\Q\otimes\E(Y,T)$ and codomain $\R$,  defined by $r\otimes\lambda\mapsto r\lambda$ for all $(r,\lambda)$ in $\Q\times\E(Y,T)$, is a $\Q$-linear monomorphism. We only need to check that the image of this map is   equal to 
  $\Bigl\{0\Bigr\}\cup\left\{t\in\R\setminus\{0\}\colon T^{{\frac{1}{t}}}\text{ is not minimal}\right\}$.
  For this purpose, we will prove that for every nonzero real number $t$, the time-$\frac{1}{t}$ map is not minimal if and only if there is $(r,\lambda)$ in $\left(\Q\times\E(Y,T)\right)\setminus\{(0,0)\}$ such that $t=r\lambda$. 
  
  Let $t$ be a nonzero real number such that $T^{\frac{1}{t}}$ is not minimal. Assume that $t$ is positive (If $t$ is negative then we work with $-t$ to obtain $-t=r\lambda$ and so $t=-r\lambda$.) Let $M$ be a minimal subset for $T^{\frac{1}{t}}$ (cf.\ \cite[Theorem~5.2]{W}). Put
    \[ 
       H=\left\{s\in\R \colon T^s\left(M\right)=M\right\}.
    \]
 
 As $H$ is a nonempty closed subgroup of $\R$ and $H\not=\R$ (because $(Y,T)$ is minimal), there is $t_0>0$ such that $H=t_0\Z$ (see e.g.\ \cite[Lemma~II.1.11]{V}). Then 
 \[ 
          \bigcup_{0\leq s<t_0} T^s(M)=T\left(M\times[0,t_0]\right)
  \]
 is a closed, nonempty and  $T$-invariant subset of $Y$.  As $T$ is minimal, we obtain that 
 \begin{equation}\label{partition}
     Y= \bigcup_{0\leq s<t_0} T^s(M).
 \end{equation}
Furthermore, for $0\leq s_1,s_2 <t_0$ with $s_1\not=s_2$ one has $T^{s_1}(M)\cap T^{s_2}(M)=\emptyset$. This together with the equality (\ref{partition}) lead us to conclude that $\{T^s(M)\colon 0\leq s < t_0\}$ forms a partition of $Y$.
Define a continuous function $\chi\colon Y\rightarrow S^1$ by the formula
 \[
     \chi(T^s(x))=e^{ 2\pi i {{\frac{1}{t_0}}}s },
 \]
for all  $x\in M$ and for all $0\leq s<t_0$. Then $\chi$ is an eigenvector of $(Y,T)$ with eigenvalue $\ts{\frac{1}{t_0}}$. Since $\ts{\frac{1}{t}}=mt_0$ for some $m\in\Z\setminus\{0\}$, we get that $t={\textstyle{ \frac{1}{m}\frac{1}{t_0} }}$, where $\left({\textstyle {\frac{1}{m},\frac{1}{t_0} }}\right)$ is in $\Q\times\E(Y,T)$.

To prove  the converse, let $(r,\lambda)$ be in $\left( \Q\times\E(Y,T)\right)\setminus\{(0,0)\}$. Assume that $\chi$ is a (nonconstant) eigenvector for $\lambda$.  Since $\chi\circ T^t =e^{2\pi i \lambda t}\chi$ for all $t\in \R$, we get that $\chi\circ T^{\frac{1}{\lambda}}=\chi$ and so $T^{\frac{1}{\lambda}}$ is not minimal, cf.\ \cite[Theorem~5.3]{W}. Since $Y$ is connected (because $(Y,T)$ is minimal, cf.\ Lemma~\ref{MinImpliesConnected}) we use Lemma~\ref{min} to conclude that $T^{\frac{1}{r\lambda}}$ is not minimal, as desired.
\end{proof}

Given a discrete flow $(X,S)$ and a strictly positive continuous function $f\colon X\to \R$, consider the function $\alpha_f\colon X\times\Z\to\R$ given by the formula
\[ 
  \alpha_f(x,n)=\begin{cases}
                        \phantom{-} \sum_{i=0}^{n-1}f(S^i(x)) & \text{if } n>0\\
                        \phantom{-}0   &\text{if } n=0\\
                        -\sum_{i=1}^{-n}f(S^{-i}(x))  &\text{if } n<0,
                        \end{cases}
\]
for all $(x,n)$ in $X\times\Z$.
Then $\R$ and $\Z$ both act on the product space $X\times\R$ by setting
\begin{eqnarray*}
       (x,s)\cdot n &=&\left(S^n(x),s-\alpha_f(x,n)\right),\text{ and }\\
       t\cdot(x,s) &=&(x,s+t),
\end{eqnarray*}
for all $(x,s)$ in $X\times\R$ and for all $(n,t)$ in $\Z\times\R$.
The two actions commute and so we let $Y_{S,f}$ be the quotient space $(X\times\R)\slash\Z$ and we let $T_{S,f}$ be the action of $\R$ on $Y_{S,f}$.
The resulting continuous flow $(Y_{S,f},T_{S,f})$ is called the {\it suspension flow  with base $S$ and ceiling function $f$}.
For each $(x,s)$ in $X\times\R$, we will denote by 
$\left[x,s\right]_f$  the image of $(x,s)$ in $Y_{S,f}$.  In this way, for the action $T_{S,f}$ of $\R$ on $Y_{S,f}$, we write $T^{t}_{S,f}\left(\left[x,s\right]_f\right)=\left[x,s+t\right]_f$ for all $t$ in $\R$ and for all $(x,s)$ in  $X\times\R$.  We suppress the subindexes $S$ and $f$ from the notation as permited by context. When $f$ is the constant function $1$, we refer to the induced flow as   the {\it standard suspension of } $S$.   A suspension flow is minimal if and only if its base is minimal. See e.g.\ \cite[Section III.5.5]{V} for details. For future reference we now prove the following.

\begin{lemma}\label{lemmamin}
  Let $(X,S)$ be a minimal dynamical system. If $t$ is a real number for which the time-$t$ map on the standard suspension of $S$ is minimal, then  $t$ is irrational.
\end{lemma}
\begin{proof}
Let $(Y,T)$ be the standard suspension flow of $(X,S)$. The time-$0$ map is the identity map and so it is not minimal. If $t={\textstyle{\frac{p}{q}}}$ is a nonzero rational, with $p$ and $q>0$ relatively prime,  then $M=\pi \left( \bigcup_{i=0}^{q-1} X\times\left\{{\textstyle{\frac{i}{q}}}\right\} \right)$ is a closed, nontrivial and $T^{t}$-invariant subset of $Y$, where $\pi\colon X\times \R\to Y$ is the canonical quotient map. Hence  the time-$t$ map on the standard suspension of $S$ is not minimal, as was to be proved. 
\end{proof}

Given a dynamical system $(X,S)$, we say that a continuous function $\xi\colon X\to S^1$ is an eigenvector for $S$ if there is $\lambda\in\R$ such that $\chi\left(S\left(x\right)\right) =e^{2\pi i \lambda}\/ \chi\left(x\right)$ for all $x\in X$. We then call $e^{2\pi i \lambda}$ the eigenvalue associated with $\chi$ or say that $\chi$ is an eigenvector with eigenvalue $e^{2\pi i \lambda}$. In analogy to the case of continuous flows, Definition~\ref{defOFeigenvalue}, we denote by $\E(X,S)$ the set of all eigenvalues for $(X,S)$. When $(X,S)$ is minimal, the set $\E(X,S)$ is actually a countable subgroup of $\T^1$, cf.\ \cite[Theorem~5.17]{W}. The following gives us a relation between the eigenvalues of a minimal dynamical system and the eigenvalues of its standard suspension flow.

\begin{lemma}\label{compareEigenvalues}
 Let $(X,S)$ be a minimal dynamical system and let $(Y,T)$ be the standard suspension of $S$. Then
        \[
           \E(X,S)=\phi(\E(Y,T)),
        \]
where $\phi\colon \R\to \T^1$ is the  canonical map $\R\ni x\mapsto e^{2\pi i x}\in\T^1$. 
\end{lemma}
\begin{proof}
   Suppose that $\lambda$ is a real number such that $e^{2\pi i\lambda}$ is an eigenvalue of $(X,S)$ with eigenvector $\chi$. Then $\lambda$ is an eigenvalue of $(Y,T)$ with eigenvector $\widetilde{\chi}\colon Y\to S^1$ defined by $\widetilde{\chi}\left(\left[x,s\right]\right)=e^{2\pi i s\lambda}\chi\left(x\right)$ for all $(x,s)\in X\times \R$.
Conversely, if $\lambda\in\R$ is an eigenvalue of $(Y,T)$ with eigenvector $\chi$ then $e^{2\pi i \lambda}$ is an eigenvalue of $(X,S)$ with eigenvector $\underline{\chi}\colon 
X\to S^1$ defined by $\underline{\chi}\left(x\right)=\chi\left(\left[x,0\right]\right)$.
\end{proof}

As mentioned in the introduction, Schwartzman \cite{Sch} introduced, for a given continuous flow $(Y,T)$ with $T$-invariant probability measure $\nu$, an asymptotic cycle $A\sb\nu \colon H^1(Y,T)\to \R$, where $H^1(Y,\Z)$ is the first \v{C}ech cohomology group with integer coefficients. We regard $H^1(Y,\Z)$ as the quotient group formed by the multiplicative group of continuous functions from $Y$ to $S^1$ modulo the ones that can  be expressed as $e^{2\pi i h\left(y\right)}$ for some continuous function $h\colon Y\to \R$. The subset $C$ consisting of all elements $[f]$ in $H^1(Y,\Z)$ such that $A\sb\nu([f])>0$, for all $T$-invariant probability measure $\nu$ on $Y$, is a (multiplicative) semigroup. Using the terminology of J.\ Packer \cite{Packer}, we will call this semigroup $C$ the {\em positive Schwartzman cone} of $H^1(Y,\Z)$. 

\begin{definition}
We define the following subset of $H^1(Y,\Z)$, which is an abelian (multiplicative) semigroup with neutral element.
\[
H^1(Y,\Z)^+=C\cup\left\{ \left[1\right]\right\}.
\]
Here $C$ is the positive Schwartzman cone of $H^1(Y,\Z)$ and $1$ is the constant function $Y\ni y \mapsto 1\in S^1$. 
\end{definition}

Let $(X,S)$ be a dynamical system and let $(Y,T)$ be the standard suspension of $S$. We will denote by 
\[
\Phi\colon K_1(C(Y))\to K_0\left(\cprod{C(X)}{S}\right) 
\] 
the Connes' Isomorphism \cite[Corollary~V.6]{C}.

 A dynamical system $(X,S)$ induces a C*-algebra $\cprod{C(X)}{\alpha}$, called the crossed product associated to $(X,S)$, which can be described as the universal C*-algebra generated by $C(X)$ and a unitary $u$ satisfying $ufu^{-1}=f\circ S^{-1}$, for all $f\in C(X)$. There is a one to one correspondence between  $S$-invariant Borel probability measures $\mu$ on $X$ and normalized traces $\tau\sb\mu$ on 
 $\cprod{C(X)}{S}$, given by 
 $\tau\sb\mu\left( \sum_kf_ku^k\right)=\int_Xf_0\,\d\mu(x)$.  A trace $\tau$ on a C*-algebra $A$ induces a natural homomorphism,  which we denote again by $\tau$, from the group $K_0(A)$ to $\R$.


\section{The Results}
We prove here our main results.
Throughout this section $(X,S)$ will denote a minimal dynamical system and $f\colon X\to \R$ will be a strictly positive, continuous function. 

Let $(Y_{S,f},T_{S,f})$ and $(Y_{S,1},T_{S,1})$ denote the suspension flows with the same base $S$ and with ceiling functions $f$ and $1$, respectively.  For each  $\lambda$ in $\E(Y_{S,f},T_{S,f})^+$ and eigenvector $\chi$ of $\lambda$, we 
consider a continuous map $U\sb\lambda\colon X\times\R\to S^1$ defined by $U\sb\lambda(x,s)=\chi\left(   \left[ x,sf(x)\right]_f \right)$ for all $(x,s)\in X\times\R$. Given $x$ in $X$ we have, by definition, that
$
   \left[ x, f(x)\right]_f = [ S(x),0]_f,
$
and so $U\sb\lambda(x,1)=U\sb\lambda(S(x),0)$. Hence the map $U\sb\lambda$ induces a continuous map with domain $Y_{S,1}$ and codomain $S^1$, which we denote again by $U\sb\lambda$,  given by the formula
\begin{align*}
U\sb\lambda \left(  \left[ x,s\right]_1 \right)
     &= \chi \left(  \left[x,sf(x)\right]_f  \right)\\
     &= \chi\left( T_{S,f}^{sf(x)} \left( \left[ x,0 \right]_f \right) \right)\\
    &=  e^{2\pi i \lambda sf(x)}\chi\left(\left[x,0\right]_f\right),
\end{align*}
for all $(x,s)\in X\times\R$. In the following lemma we show that the class $[U\sb\lambda]$  in $H^1(Y_{S,1},\Z)$ does not depend on the choice of the eigenvector $\chi$ of $\lambda$; moreover, $[U\sb\lambda]$  belongs to the positive Schwartzman cone of $H^1(Y_{S,1},\Z)$,  except when $\lambda$ is zero.

\begin{lemma}\label{gamma}
Let $(Y_{S,f},T_{S,f})$ and $(Y_{S,1},T_{S,1})$ denote the suspension flows with the same base $S$ and with ceiling functions $f$ and $1$, respectively.  For each $\lambda$ in $\E(Y_{S,f},T_{S,f})^+$ and eigenvector $\chi$ of $\lambda$,  consider the map $U\sb\lambda\colon Y_{S,1}\to S^1$ as defined above. Then the class $[U\sb\lambda]$  in $H^1(Y_{S,1},\Z)$ does not depend on the choice of the eigenvector $\chi$ of $\lambda$; moreover, $[U\sb\lambda]$  belongs to the positive Schwartzman cone of $H^1(Y_{S,1},\Z)$,  except when $\lambda$ is zero.  
\end{lemma}
\begin{proof}
Suppose that $\chi_1$ and $\chi_2$ are two different eigenvectors for $\lambda$. Using Proposition~\ref{eigen}, there is a constant $c$ in $S^1$ such that $\chi_1=c\chi_2$.
Therefore, for each $(x,s)$ in $X\times\R$, it follows that 
$\chi_1\left(\left[x,sf(x)\right]_f\right)=c\/ \chi_{2}\left(\left[x,sf(x)\right]_f\right)$. This proves that the class of $U\sb\lambda$ in 
$H^1(Y_{S,1},\Z)$ does not depend on the choice of the eigenvector $\chi$ of $\lambda$, as wanted.   

 Now, observe that for each $(x,s)$ in $X\times\R$ we can write
 \[
 U\sb\lambda(\left[x,s\right]_1)=\underline{\chi}(x)e^{2\pi i sh(x)},
 \] 
 where 
 $\underline{\chi}\colon X\to S^1$ and $h\colon X\to \R$ 
 are two continuous functions on $X$ which are defined by the formulas
 $\underline{\chi}(x)=\chi\left(\left[x,0\right]_f\right)$ 
 and 
 $h(x)=\lambda f(x)$, for all $x\in X$. Furthermore, for each $x$ in $X$, the functions $\underline{\chi}$ and $h$ satisfy
 \[
    \underline{\chi}\left(S(x)\right)=\underline{\chi}\left(x\right) e^{2\pi i h(x)}.
 \]
 Using \cite[Theorem~1.1]{Packer2}, we have that $[U\sb\lambda]$ is in the positive Schwartzman cone of $H^1(Y,\Z)$ if and only if $\int_Xh(x)\,\d\mu(x)$ is strictly positive for every $S$-invariant Borel probability measure $\mu$ on $X$.
We  claim that this is the case, unless $\lambda=0$. Since $f$ is strictly positive, it follows that $\int_Xf(x)\,\d\mu(s)$ is strictly positive for every $S$-invariant Borel probability measure $\mu$ on $X$. Therefore the integral $\int_X h(x)\,\d\mu(x)=\int_X\lambda f(x)\,\d\mu(x)$ is strictly positive for every $S$-invariant Borel probability measure $\mu$ on $X$, unless $\lambda$  is zero. This proves the claim and completes the proof of the lemma.
\end{proof}

\begin{proposition}\label{semigroup}
Let the notation be as in Lemma~\ref{gamma}. Define a function 
\[
\Gamma_{S,f}^+\colon \E(Y_{S,f},T_{S,f})^+\to H^1(Y_{S,1},\Z)^+
\]
by $\Gamma_{S,f}^+(\lambda)=[U\sb\lambda]$ for all $\lambda\in\E(Y_{S,f},T_{S,f} )^+$. Then the map $\Gamma_{S,f}^+$ is a homomorphism of abelian semigroups with neutral elements.
\end{proposition}
\begin{proof}
    
   Let $\lambda_1$ and $\lambda_2$ be two elements in $\E(Y_{S,f},T_{S,f})^+$ with eigenvectors $\chi_1$ and $\chi_2$, respectively.  Then $\chi_1\chi_2$ is an eigenvector for $\lambda_1+\lambda_2$. Since by Lemma~\ref{gamma} the class of $\left[ U_{\lambda_1+\lambda_2}\right]$ does not depend on the choice of an eigenvector for $\lambda_1+\lambda_2$, we get
 \begin{align*}
     \Gamma_{S,f}^+\left(\lambda_1+\lambda_2\right) &= \left[U_{\lambda_1+\lambda_2}\right]\\
        &= \left[\left(\chi_1\chi_2\right)\left(\left[x,sf(x)\right]_f\right)\right]\\
        &= \left[\chi_1\left(\left[x,sf(x)\right]_f\right)\right]\, \left[\chi_2\left(\left[x,sf(x)\right]_f\right)\right]\\
        &=\left[U_{\lambda_1}\right]\, \left[U_{\lambda_2}\right]\\
        &=\Gamma_{S,f}^{+}\left(\lambda_1\right)\ \Gamma_{S,f}^+\left(\lambda_2\right).
 \end{align*} 
    Hence $\Gamma_{S,f}^+$ is a  homomorphism, as wanted. 
\end{proof}

Observe that, given a continuous flow $(Y,T)$, the group $\E\left(Y,T\right)$ can be seen as the Grothendieck group of $\E\left(Y,T\right)^+$. Furthermore, the Grothendieck group of (the multiplicative semigroup) $H^1(Y,\Z)^+$ is a subgroup of $H^1(Y,\Z)$. In particular, we can  consider the induced group homomorphism $\Gamma_{S,f}$ of the map $\Gamma_{S,f}^+$ of Proposition~\ref{semigroup}. We state this observation in the following.

\begin{proposition}\label{group}
Let the notation be as in Proposition~\ref{semigroup}. There exists a unique group homomorphism $\Gamma_{S,f}\colon\E(Y_{S,f},T_{S,f})\to H^1(Y_{S,1},\Z)$ that makes the following diagram commute. \vspace{-.2cm}
\[ 
  \xymatrix{ %
     \E\left(Y_{S,f},T_{S,f}\right)^+\ar[rr]^(.55){\Gamma_{S,f}^+}\ar[d]_{\iota} &&   
                       H^1\left(Y_{S,1},\Z\right)^+     \ar[d]^{\iota}\\
     \E\left(Y_{S,f},T_{S,f}\right)\ar@{-->}[rr]^(.55){\Gamma_{S,f}}  && H^1\left(Y_{S,1},\Z\right)  \\
                   }
\]
\end{proposition}
\begin{proof}
Let $G$ be the Gro\-then\-dieck group of $H^1\left(Y_{S,1},\Z\right)^+$. Then the existence of a unique homomorphism $\Gamma_{S,f}\colon  \E(Y_{S,f},T_{S,f})\to G\subset H^1(Y_{S,1},\Z)$ such that the diagram above is commute is a consequence of the fact that $\E\left(Y_{S,f},T_{S,f}\right)$ is the Grothendieck group of $\E\left( Y_{S,f},T_{S,f}\right)^+$, see e.g.\ a corollary in \cite[Appendix~D]{W-O}, p.297. 
\end{proof}

Considering $H^1(Y_{S,1},\Z)$ as  a direct summand of $K_1(C(Y_{S,1}))$, we can now bring $K$-theory into the picture by composing the map $\Gamma\colon\E(Y_{S,f},T_{S,f})\to H^1(Y_{S,1},\Z)$ in Proposition~\ref{group} with the Connes' Isomorphism $\Phi\colon K_1(C(Y_{S,1}))\to K_0\left(\cprod{C(X)}{S}\right)$. Relying on the work of Packer \cite{Packer2}, we will show next  how the composition map
$\Phi\circ\Gamma\colon \E\left(Y_{S,f},T_{S,f}\right)\to K_0\left(\cprod{C(X)}{S}\right)$ has a left inverse given by a map induced by a trace on $\cprod{C(X)}{S}$.

\begin{proposition}\label{taumu}
   Let the notation be as in Proposition~\ref{group}. If $\mu$ is an $S$-invariant Borel probability measure on $X$ then the map
   \[
       \tau_{\mu,f\colon}K_0\left(\cprod{C(X)}{S}\right)\to \R,
   \]
defined by $\tau_{\mu,f}(x)=\ts{\frac{\tau\sb\mu\left(x\right)}{\tau\sb\mu\left(f\right)}}$ for all $x$ in 
$K_0\left(\cprod{C(X)}{S}\right)$, is a left inverse for the map $\Phi\circ\Gamma_{S,f}\colon \E\left(Y_{S,f},T_{S,f}\right)\to K_0\left(\cprod{C(X)}{S}\right)$.
\end{proposition}
\begin{proof}
We have to show that $\left(\tau_{\mu,f}\circ\Phi\circ\Gamma_{S,f}\right)(\lambda)=\lambda$ for all  $\lambda$ in $\E(Y_{S,f},T_{S,f})$. Since the maps involved are group homomorphisms, it will suffice to do it for strictly positive $\lambda$. Let $\lambda\in\E\left(Y_{S,f},T_{S,f}\right)^+\setminus\{0\}$ and say $\chi$ is an eigenvector for $\lambda$. By Lemma~\ref{gamma} and Proposition~\ref{semigroup},  we have that $\Gamma_{S,f}(\lambda)=[U\sb\lambda]=\left[\chi\left( \left[x,0\right]_f  \right)e^{2\pi i \lambda sf(x)}\right]$ is in the positive Schwartzman cone of $H^1(Y_{S,1},\Z)$. Then, by \cite[Theorem~3.4]{Packer2}, there is a projection $p$ in $M_n\left(\cprod{C(X)}{S}\right)$ such that $\Phi([U\sb\lambda])=[p]$. Furthermore, \cite[Theorem~2.3]{Packer2} gives us that 
\[
\tau\sb\mu(p)=\int_X\lambda f(x)\,\d\mu(x)=\lambda\tau\sb\mu(f).
\]
 Thus
\[
\left(\tau_{\mu,f}\circ\Phi\circ\Gamma_{S,f}\right)\left(\lambda\right)=\tau_{\mu,f}
      \left(\Phi\left(\left[U\sb\lambda\right]\right)\right)
    =\tau_{\mu,f}\left(p\right)  
    =\frac{\tau\sb\mu\left(p\right)}{\tau\sb\mu\left(f\right)}
    =\lambda.
\]
\end{proof}

We are ready to prove our main two results.

\begin{theorem}\label{step2}
Let $(Y,T)$ be a minimal suspension flow with base  $(X,S)$ and with (strictly positive, continuous) ceiling function $f\colon X\to \R$. Let $(Y_S,T_S)$ denote the standard suspension of $S$. Consider the homomorphism $\Gamma\colon \E(Y,T)\to H^1(Y_S,\Z)$ of Proposition~\ref{group} and put
  \[
      \E K(X,S,f)=\left(\Phi\circ \Gamma\right) \left(  \E(Y,T)  \right).
  \]
If $\mu$ is an  $S$-invariant Borel probability measure on $X$ then the 
map
  \[
    \tau_{\mu,f}\colon K_0\left(\cprod{C(X)}{S}\right)\to\R
  \]
defined by the formula $\tau_{\mu,f}(x)={\textstyle \frac{\tau\sb\mu(x)}{\tau\sb\mu(f)}}$, for all $x\in K_0\left(\cprod{C(X)}{S}\right)$, satisfy $\tau_{\mu,f}(\E K(X,S,f))=\E(Y,T)$. 
Moreover, the map
\[\tau_{\mu,f}\colon \E K(X,S,f)\to \E(Y,T)\] is an order isomorphism of ordered abelian groups, where $\E(Y,T)$ inherits the order from $\R$.
\end{theorem}
\begin{proof}   
 Using Proposition~\ref{taumu}, we get  that
\begin{equation}\label{isom}
       \Phi\circ\Gamma\colon \E\left(Y,T\right)\to\E K\left(X,S,f\right)\subset K_0\left(\cprod{C(X)}{S}\right) 
\end{equation}
is an isomorphism with inverse $\tau_{\mu,f}\colon \E K\left(X,S,f\right)\to \E\left(Y,T\right)$.
In fact, the isomorphism in (\ref{isom}) can be seen to be an order isomorphism since, by  \cite[Theorem~3.4]{Packer2} (see also \cite[Remark~3.6]{Packer2}), we have that
$\left(\Phi\circ\Gamma\right)(\E (Y,T)^+)$ lies in the positive cone $K_0\left(\cprod{C(X)}{S}\right)^+$ of $K_0\left(\cprod{C(X)}{S}\right)$. Thus $\tau_{\mu,f}\colon \E K(X,S,f)\to \E(Y,T)$ is also an order isomorphism, as was to be proved.
\end{proof}

\begin{theorem}\label{withCTI}
Let $(Y,T)$ be a minimal suspension flow with base  $(X,S)$ and with (strictly positive, continuous) ceiling function $f\colon X\to \R$. Let $\mu$ be an $S$-invariant Borel probability measure on $X$. Consider the map with domain $\Q\otimes\E K(X,S,f)$ and codomain $\R$ defined by $r\otimes x \mapsto r{\ts {\frac{\tau\sb\mu(x)}{\tau\sb\mu(f)}}}$ for all $(r,x)$ in $\Q\times \E K (X,S,f)$, where $\E K(X,S,f)$ is as in Theorem~\ref{step2}. This map is a $\Q$-linear monomorphism with range equal to
\[ 
\Bigl\{0\Bigr\}\cup\left\{t\in\R\setminus\{0\}\colon T^{\frac{1}{t}}\text{ is not minimal}\right\}.
\]
Hence the set above is a countable $\Q$-linear subspace of $\R$ isomorphic to
$ \Q\otimes\E K(X,S,f)$.
\end{theorem}
\begin{proof}
    Since  Theorem~\ref{step2} says that $\E K(X,S,f)$ is isomorphic to $\E(Y,T)$ via the map defined by $x\mapsto { \textstyle\frac{\tau\sb\mu\left(x\right)}{\tau\sb\mu\left(f\right)} }$ for all $x\in \E K(X,S,f)$, we obtain that the map defined by 
    \[r\otimes x\mapsto r\otimes\textstyle\frac{\tau\sb\mu\left(x\right)}{\tau\sb\mu\left(f\right)},
    \] 
for all $(r,x)$ in $\Q\times \E K(X,S,f)$, gives an isomorphism of $\Q\otimes \E K(X,S,f)$ with $\Q\otimes \E(Y,T)$. On the other hand, Proposition~\ref{step1} asserts that 
$\Q\otimes \E(Y,T)$ is isomorphic to the $\Q$-linear space
 $
   \Bigl\{0\Bigr\}\cup\left\{t\in\R\setminus\{0\}\colon T^{\frac{1}{t}}\text{ is not minimal}\right\}
 $ 
via the map defined by \[r\otimes \lambda \mapsto r\lambda\] for all $(r,\lambda)\in \Q\times \E(Y,T)$. Combining these two isomorphism we get the desired result.
\end{proof}

When we specialize to the case for which $f$ is the constant function $1$ (that is, when we consider the standard suspension of $S$), we can say the following.

\begin{corollary}\label{mainCor}
Let $(X,S)$ be a minimal dynamical system. If $t$ is a nonzero real number such  that the time-$\textstyle\frac{1}{t}$ map on the standard suspension of $S$ is not minimal then there exists $(r_1,r_2)\in \Q\times\Q$ and there exists a Rieffel projection $p$ in $\cprod{C(X)}{S}$ such that for all $S$-invariant Borel probability measure $\mu$ on $X$ we have   $t=r_1+r_2\tau\sb\mu\left(p\right)$. 
\end{corollary}
\begin{proof}
Denote by $(Y,T)$ the standard suspension flow of $(X,S)$.
Since $T^{\frac{1}{t}}$ is not minimal, Proposition~\ref{step1} implies the existence of  $(r,\lambda)$ in $\Q\times\E(Y,T)$ such that $t=r \lambda$. Let $n$ and $\gamma$ denote the integral and fractional parts of $\lambda$, respectively. To avoid trivialities, assume that $\lambda$ is irrational (else we may choose $p$ to be trivial). Then $\gamma$ lies in the open interval $(0,1)$. Observe that  $\Z\subset \E(Y,T)$, cf.\ Lemma~\ref{compareEigenvalues}. Therefore, as $\E(Y,T)$ is a subgroup of $\R$ (cf.\ Proposition~\ref{eigen}), the fractional part $\gamma$ of $\lambda$ belongs to $\E(Y,T)^+\setminus\{0\}$. By Theorem~\ref{step2}, there exists a projection
$p$ in $M_n\left(\cprod{C(X)}{S}\right)$ (in fact $[p]=\Phi([U\sb\gamma])$) such that  $\gamma=\tau\sb\mu(p)=\tau_{\mu,1}(p)$ for all $S$-invariant Borel probability measure $\mu$ on $X$. Combining now Theorems 3.2 and 3.4 in \cite{Packer2}, we conclude that $p$ is a Rieffel projection. Thus $t=rn+r\tau\sb\mu\left(p\right)$, as desired.
\end{proof}

\begin{corollary}
Let $(X,S)$ be a minimal Cantor system. If $t$ is a nonzero real number such that the time-$\textstyle\frac{1}{t}$ map on the standard suspension of $S$ is not  minimal then there exists  $(r_1,r_2)\in\Q\times\Q$ and there exits an open and closed subset $F$ of $X$ such that  for all $S$-invariant Borel probability measure $\mu$ on $X$ we have $t=r_1+r_2\mu\left(F\right)$.
\end{corollary}
\begin{proof}
By Corollary~\ref{mainCor}, there is $(r_1,r_2)\in\Q\times\Q$ and there is a Rieffel projection $p$ in $\cprod{C(X)}{S}$ such that for all $S$-invariant Borel probability measure $\mu$ on $X$ we have $t=r_1+r_2\tau\sb\mu(p)$. To avoid trivialities we assume $t$ is irrational (else we may choose $F$ to be  trivial). Therefore $0 < \tau\sb\mu(p)<1$ for all $S$-invariant Borel probability measure $\mu$ on $X$. Furthermore,  since $p$ is a projection in $\cprod{C(X)}{S}$, then $p$ can be regarded as a continuous function in $C(X,\Z)$, cf. \cite[Theorem~1.1]{Putnam}. Using \cite[Lemma~2.4]{GW}, we find an open and closed subset $F$ of $X$ such that $\tau\sb\mu(p)=\mu(F)$ for all $S$-invariant Borel probability measure $\mu$ on $X$. Thus $t=r_1+r_2\mu(F)$ for all  $S$-invariant Borel probability measure $\mu$ on $X$, as wanted. 
\end{proof}

Extending a result of Rieffel, Pimsner and Voiculescu \cite{R}, N.\ Riedel \cite{Riedel} showed that the range of the unique trace on the $K_0$-group of a minimal group rotation can be described in terms of the eigenvalues  of this group rotation. This result was later proved by other authors using different techniques (see e.g.\ \cite[Theorem~2.5]{Packer} and \cite[Theorem~IX.11]{E}). For general minimal dynamical systems we get the following.

\begin{corollary}\label{Riedel}
   Let $(X,S)$ be a minimal dynamical system and let $\mu$ be an $S$-invariant Borel probability measure on $X$. Then
   \[
      \E(X,S)\subset\phi\left( \tau\sb\mu\left( K_0\left(\cprod{C(X)}{S} \right) \right)\right),
   \]
where $\phi\colon\R\to\T^1$ is the canonical map $\R\ni x\mapsto e^{2\pi i x}\in\T^1.$   
\end{corollary}
\begin{proof}
Let $(Y,T)$ be the standard suspension of $S$. Then
  \begin{align*}
      \E(X,S) &= \phi (\E(Y,T)) \text{ by Lemma~\ref{compareEigenvalues}} \\
     &\subset \phi(\tau_{\mu,1}( K_0\left(\cprod{C(X)}{S}\right))) \text{ by Theorem~\ref{step2}}\\
     &= \phi\left( \tau\sb\mu\left( K_0\left(\cprod{C(X)}{S}\right) \right) \right)
         \text{ by definition of $\tau_{\mu,1}$}.
  \end{align*}
\end{proof}

 We show next that the $K$-theory of a minimal Cantor system keeps track of  the eigenvalues of every dynamical system in its orbit equivalence class (cf.\ \cite[Definition~1.2]{GPS}).

\begin{corollary}\label{eigenvaluesInTrace}
Let $(X,S)$ be a minimal Cantor system and let $\mu$ be an $S$-invariant Borel probability measure on $X$. If $(X_1,S_1)$ is orbit equivalent to $(X,S)$ then
\[
    \E(X_1,S_1)\subset\phi \left(\tau\sb\mu\left( K_0\left(\cprod{C(X)}{S}\right)\right) \right),
\]
where $\phi\colon\R\to\T^1$ denotes the canonical map $\R\ni x\mapsto e^{2 \pi i x}\in\T^1$.
\end{corollary}
\begin{proof}
Since $(X_1,S_1)$ is orbit equivalent to  $(X,S)$, there is a homeomorphism $F\colon X_1\to X$ carrying the $S_1$-invariant Borel probability measures on $X_1$ onto the $S$-invariant Borel probability measures on $X$. Furthermore, this homeomorphism induces an order isomorphism of $K_0\left(\cprod{C(X_1)}{S_1}\right)\slash\ \mathrm{Inf}\left(K_0\left(\cprod{C(X_1)}{S_1} \right)\right)$ with $K_0\left(\cprod{C(X)}{S}\right)\slash\ \mathrm{Inf}\left(K_0\left(\cprod{C(X)}{S} \right)\right)$, cf.\ 
\cite[Theorem~2.2]{GPS}. Let $\nu$ be the $S_1$-invariant measure corresponding to
$\mu$ under $F$. Then
\begin{align*}
\tau\sb\nu\left(K_0\left(\cprod{C(X_1)}{S_1}\right)\right)
   &=  \tau\sb\nu \bigl(K_0\left(\cprod{C(X_1)}{S_1}\right)\slash\ \mathrm{Inf}\left(K_0\left(\cprod{C(X_1)}{S_1} \right)\right) \bigr)\\
   &= \tau\sb\mu\bigl(
  K_0\left(\cprod{C(X)}{S}\right)\slash\ 
   \mathrm{Inf}\left(K_0\left(\cprod{C(X)}{S} \right)\right)\bigr)\\
   &= \tau\sb\mu\left(
  K_0\left(\cprod{C(X)}{S}\right)
                         \right)
\end{align*}
Combining the above equality with Corollary~\ref{eigenvaluesInTrace} we obtain
\begin{align*}
  \E (X_1,S_1) &\subset \phi\left(  \tau\sb\nu \left(  K_0(\cprod{C(X_1)}{S_1}) \right) \right)\\
       &= \phi\left(  \tau\sb\mu \left(   K_0 ( \cprod{C(X)}{S} )\right)\right),
\end{align*}
as desired.
\end{proof}


\section{Examples}
In this section we include some examples illustrating our results. We start by showing dynamical systems for which  every irrational time map on its standard suspensions is minimal. To contrast these, we also show some dynamical system for which some irrational time map on its standard suspension is not minimal. Our next examples will present dynamical systems $(X,S)$ one which has group $\E K(X,S,1)$ (cf.\ Theorem~\ref{step2}) equal to $K_0\left(\cprod{C(X)}{S}\right)$ and one for which $\E K(X,S,1)$ is different from $K_0\left(\cprod{C(X)}{S}\right)$. We will also give an example of two strong orbit equivalent minimal Cantor systems which have different eigenvalues. Our last example exhibits a minimal suspension flow $(Y,T)$ satisfying the following: the map $T^t$  is minimal for every nonzero real number $t$.

\begin{proposition}\label{forQ}
Let $(X,S)$ be a minimal  dynamical system. Suppose there is an $S$-invariant Borel probability measure $\mu$ on $X$ such that $\tau\sb\mu\left(K_0\left(\cprod{C(X)}{S}\right)\right) $ is contained in $\Q$. Then, for each real number $t$, the time-$t$ map on the standard suspension of $S$  is minimal if and only if  $t$ is irrational.
\end{proposition}
\begin{proof}
 Suppose  that $t$ is a real number such that the time-$t$ map on the standard suspension of $(X,S)$ is not minimal.  We show that $t$ is rational. If $t$ is zero there is nothing to prove. If $t$ is nonzero then, by Theorem~\ref{withCTI}, there is $(r,x)$ in $\Q\times\E K(Y,T,1)$ such that  $\frac{1}{t}=r\tau\sb\mu(x)$. As $\E K(X,S,1)$ is a subgroup of $K_0\left(\cprod{C(X)}{S}\right)$ and $\tau\sb\mu\left(K_0(\cprod{C(X)}{S}\right) $ is contained in $\Q$, we conclude that $\tau\sb\mu(x)$ is rational and so $t$ is rational, as wanted. The converse is Lemma~\ref{lemmamin}.
\end{proof}

\begin{example}
The rotation $R_t\colon z\mapsto ze^{2\pi i t}$ on $S^1$ is minimal if and only if $t$ is irrational.
\end{example}
This is a well know fact, see e.g.\ \cite[Proposition~III.1.4]{V}. We use our results for an alternative proof. As the unit circle $S^1$ can be regarded as the standard suspension of the one-point dynamical system, the rotation $R_t$ then corresponds to the time-$t$ map of this suspension. Since the range of the trace on the $K_0$-group of the one-point dynamical system is equal to $\Z$, we apply Proposition~\ref{forQ} to get that $R_t$ is minimal if and only if $t$ is irrational.\hfill $\Box$\medskip

Let $\{n_i\}_{i=1}^{\infty}$ be a sequence of integers, each greater than $2$. Let $X=\Pi_{i=1}^{\infty}\{0,1,\ldots,n_{i-1}\}$. Recall that an {\em odometer system} consists of the space $X$ together with the homomorphism $S$ representing addition of $(1,0,0,\ldots,0)$ with carry over to the right.

\begin{example}\label{Odometer}
Let $t$ be a real number. The time-$t$ map on the suspension of an odometer is minimal if and only if  $t$ is irrational. 
\end{example}
The image of the trace on the $K_0$-group of an odometer is contained in $\Q$, cf.\ \cite[p.332]{Putnam}. Therefore Proposition~\ref{forQ} applies to get the desired result. \hfill $\Box$
\bigskip

Given two dynamical systems $(X_1,S_1)$ and $(X_2,S_2)$, we say that  $(X_1,S_1)$ is an extension of $(X_2,S_2)$ if there exists a continuous surjective map $F\colon X_1\rightarrow X_2$ such that $F\circ S_1=S_2\circ F$. In such case we call $F$ an extension map. We say that $(X_1,S_1)$ is an almost one to one extension of $(X_2,S_2)$ if there is an extension map $F\colon X_1\rightarrow X_2$ such that the set $A=\left\{ x\in X_1\colon F^{-1}\left(F(x)\right) = \{x\}\right\}$ is dense in $X_1$. We call such $F$ an almost one to one extension map.  It is easy to check that when $(X_1,S_1)$ is an almost one to one extension of $(X_2,S_2)$ \ then $S_1$ is minimal if and only if $S_2$ is minimal. 

\begin{lemma}\label{almost1-1}
If $(X_1,S_1)$ is an almost one to one extension of $(X_2,S_2)$ then, for each real number $t$, 
the time-$t$ map on the suspension of $(X_1,S_1)$ is minimal if and only if the time-$t$ map on the suspension of $(X_2,S_2)$ is minimal.
\end{lemma}
\begin{proof}
   Let $(Y_1,T_{1})$ and $(Y_2,T_2)$ be the standard suspensions of $S_1$ and $S_2$, respectively.   
   The lemma will follow once we prove that, for each $t\in\R$,  the dynamical system $(Y_1,T_{1}^{t})$ is an almost one to one extension of $(Y_2,T_{2}^{t})$.    
   Let  $F\colon X_1\rightarrow X_2$ be an almost one to one extension map. Consider the continuous surjective map with domain $X_1\times \R$ and codomain $Y_2$ given by $(x,s)\mapsto [F(x),s]$, for all $(x,s)$ in $X_1\times\R$. Since for each $x$ in $X_1$ we have that $[F(x),1]=[S_2(F(x)),0]=[F(S_1(x)),0]$, this map induces a continuous surjective map, which we denote $\widetilde{F}=\widetilde{F_t}\colon Y_1\rightarrow Y_2$, given by the formula $\widetilde{F}([x,s]) = [F (x),s]$, for all $(x,s)$ in $X_1\times\R$.  Since for each $(x,s)$ in $X_1\times\R$ we have
 $ T_{2}^{t}\circ\widetilde{F}\left([x,s]\right)= [F(x),s+t]= \widetilde{F}\circ T_{1}^{t}\left([x,s]\right)$, we conclude that $\widetilde{F}$ is  an extension map. Moreover, if we denote $A=\left\{ x\in X_1\colon F^{-1}\left(F(x)\right)=\{x\}\right\}$, we have that 
\begin{equation}\label{dense}
 \left\{ y\in Y_1\colon {\widetilde{F}}^{-1}\left(\widetilde{F}(y)\right)=\left\{y\right\}\right\} 
          =  \Bigl\{ [x,s]\in Y_1\colon x\in A,\, s\in\R\Bigr\}.
\end{equation}
Since $A$ is dense in $X_1$ then $A\times\R$ is dense in $X_1\times\R$ and so the right hand side in (\ref{dense}) is dense in $Y_1$. Hence the left hand side in (\ref{dense}) is  dense $Y_1$.  This proves that $\widetilde{F}$ is an almost one to one extension map from $(Y_1,T_{1}^{t})$ onto $(Y_2,T_{2}^{t})$, as was to be proved.
\end{proof}

\begin{example}\label{Toeplitz}
Let $t$ be a real number. The time-$t$ map on the standard suspension of a Toeplitz flow is minimal if and only if $t$ is irrational.
\end{example}
Since a Toeplitz flow is an almost one to one extension of an odometer (cf.\ \cite[Theorem~4]{DL}),  we combine Lemma~\ref{almost1-1} with Example~\ref{Odometer} to obtain the desired result.\hfill$\Box$\bigskip

Example~\ref{Toeplitz} together with Theorem~\ref{step1} give us that Toeplitz flows can not have irrational eigenvalues. This is in contrast with the measure-theoretical case: Toeplitz flows may have measure-theoretical irrational eigenvalues, cf. \cite{Iwanik}.
It is an interesting problem to determine when a measure-theoretical eigenvalue of a dynamical system is a ``continuous" eigenvalue (in our sense), cf.\ \cite{CDHM}.

In the following two examples we present dynamical systems for which some irrational time maps on  its standard suspensions are not  minimal.

\begin{example}\label{susp_of_Rs}
 Let $s$ be an irrational number and  let $R_s\colon z\mapsto ze^{2\pi i s}$ be the rotation by $s$ on the unit circle. For each real number $t$, the time-$\frac{1}{t}$ map on the standard suspension of $R_s$ is minimal if and only if  $t$ does not belong to $\Q+s\Q$.
\end{example}
Let $(Y,T)$ denote the standard suspension flow of $R_s$. For each $t\in\R$, we claim that $(Y,T^t)$ is conjugate to $(S^1\times S^1, R_{st}\times R_{t})$, where $R_{st}\times R_{t}$ denotes the homeomorphism of $S^1\times S^1$ given by $\left(R_{st}\times R_{t}\right) \left(z_1,z_2\right)= \left(e^{2\pi i st } z_1, e^{2\pi i t}z_2\right)$, for all $(z_1,z_2)$ in $S^1\times S^1$. Indeed, consider the map with domain $S^1\times\R$ and codomain $S^1\times S^1$ defined by $(z,r)\mapsto \left(e^{2\pi i sr}z,e^{2\pi i r}\right)$, for all $(z,r)$ in $S^1\times\R$. Since for each $z\in S^1$ we have that
$\left(e^{2\pi i s}z,e^{2\pi i }\right)=\left(R_s(z),e^0\right)$, this map induces a homeomorphism, which we denote $F=F_t\colon Y\rightarrow S^1\times S^1$, 
given by the formula $F([z,r])=\left(e^{2\pi i sr}z,e^{2\pi i r}\right)$ for all $(z,r)$ in $S^1\times\R$. We check that $F\circ T^t = \left(R_{st}\times R_t\right)\circ F$. For this purpose, let $(z,r)$ be an element in $S^1\times\R$. Then 
\begin{align*}
\left(F\circ T^t\right) ([z,r]) &= F\left([z,r+t]\right)\\
                             &= \left(e^{2\pi i s(r+t)}z, e^{2\pi i (r+t)}\right)\\
                             &= \left(R_{st}\times R_t\right)\left( e^{2\pi i sr}z,e^{2\pi i r}\right)\\
                             &= \left(\left(R_{st}\times R_t\right)\circ F\right) \left([z,u]\right).
\end{align*}
This proves the claim. Now, since $(S^1\times S^1, R_{st}\times R_{t})$ is minimal if and only if $1$, $st$ and $t$ are linearly independent over $\Q$ (see, e.g.\ \cite[Proposition~1.4.1]{K}). Then, by conjugacy, $(Y,T^{\frac{1}{t}})$ is minimal if and only $1$, $\textstyle\frac{s}{t}$ and $\textstyle\frac{1}{t}$ are linearly independent over $\Q$. This happens if and only if $t$ is not in $\Q+s\Q$, as was to be proved.
\hfill$\Box$\medskip

We now give a brief description of Denjoy systems. Recall that an orientation preserving homeomorphism of the unit circle $S^1$ without periodic points necessarily has an irrational rotation number, cf.\ \cite[Proposition 11.1.4]{K}.  A Denjoy homeomorphims  is an orientation preserving homeomorphims of $S^1$ with irrational rotation number and which is not conjugate to an irrational rotation of $S^1$. Denjoy homeomorphisms are uniquely ergodic, cf.\ \cite[Theorems 11.2.7 and 11.2.9]{K}.  When we restrict a Denjoy homeomorphism $S$ to the support $X$ of its unique invariant measure,  we obtain a strictly ergodic Cantor system $(X,S)$ which we call {\em Denjoy system}. It turns out that a Denjoy system $(X,S)$ with irrational rotation number $s$ is an almost one to one extension of the irrational rotation $R_s\colon z\mapsto e^{2\pi i s}z$ of $S^1$. Furthermore, the conjugacy class of a Denjoy system $(X,S)$ is determined by the rotation number $s$ of $S$ and a countable subset of the circle, denoted by $Q(S)$, which is invariant under rotation by an angle $2\pi i s$.   
  We refer to \cite{PSS} for a detailed account.

\begin{example}\label{Denjoy}
 Let $t$ be a real number. The time-$\frac{1}{t}$ map on the standard suspension of a Denjoy system with irrational rotation number $s$ is minimal if and only if $t$ does not belong to $\Q+s\Q$.
\end{example}

Let  $(X,S)$ be a Denjoy system with irrational rotation number $s$. Since $(X,S)$  is an almost one to one extension of the irrational rotation $R_s\colon z\mapsto e^{2\pi i s}z$ on $S^1$ then, by Lemma~\ref{almost1-1}, the time-$\textstyle\frac{1}{t}$ map on the standard suspension of $S$
 is minimal if and only if the time-$\textstyle\frac{1}{t}$ map on the standard suspension of  $R_s$ is minimal.  Using Example~\ref{susp_of_Rs}, this happens if and only if $t$ is not in $\Q+s\Q$.\hfill$\Box$

\medskip

Our next two examples will show that $K_0$-group of a minimal dynamical systems might contain more information than just the eigenvalues of its standard suspension, 
cf.\ Corollary~\ref{Riedel}.

\begin{example}\label{containedStrictly}
There are minimal dynamical systems $(X,S)$ with  $\E K(X,S,1)$ strictly contained in $K_0\left(\cprod{C(X)}{S}\right)$ 
\end{example}
Let $s_1$ and $s_2$ be real number such  that $1$, $s_1$ and $s_2$ are linearly independent over $\Q$. Consider the Denjoy systems $(X,S)$  with the following conjugacy invariants: the irrational rotation number of (X,S) is $s_1$ and 
the set $Q(S)$ is $\{ e^{2\pi i ns_1}\colon n\in\Z\}\cup \{ e^{2\pi i \left(s_2+ns_1\right)}\colon n\in\Z\}$.   By \cite[Theorem~5.3]{PSS}, the (unique) normalized trace $\tau$ on $\cprod{C(X)}{S}$ induces an (order) isomorphism between $K_0\left(\cprod{C(X)}{S}\right)$ and $\Z+s_1\Z+s_s\Z$. Let $x$ be the element in $K_0\left(\cprod{C(X)}{S}\right)$ for which $\tau(x)=s_2$. Then $x$ does not belong to $\E K(X,S,1)$; for if it does, then Proposition~\ref{step2} says that $\tau(x)=s_2$ is an eigenvalue for the standard suspension of $(X,S)$ and so, by Proposition~\ref{step1},  the time-$\textstyle\frac{1}{s_2}$  map on that suspension is not minimal. Using  Example~\ref{Denjoy} we conclude that $s_2$ belongs to $\Q+s_1\Q$, contradicting the 
fact that $1$, $s_1$ and $s_2$ are linearly independent over $\Q$.\hfill$\Box$

\begin{example}
There are minimal dynamical systems $(X,S)$ with $\E K(X,S,1)$ equal to $K_0\left(\cprod{C(X)}{S}\right)$.
\end{example}
Let $s$ be an irrational number and consider the irrational rotation $R_s\colon z\mapsto ze^{2\pi i s}$ on $S^1$.  The induced crossed product $\cprod{C(S^1)}{R_s}$ is called an {\it irrational rotation algebra} \cite{R}. Let  $(Y,T)$ denote the standard suspension flow of $\left(S^1,R_s\right)$. It is known that $\E\left(Y,T\right)$ is equal to the range of the (unique) trace on $K_0\left(\cprod{C(S^1)}{R_s}\right)$, see e.g.\ \cite[Example~1 in pg.155]{Riedel}, \cite[Theorem~2.5]{Packer} or \cite[Example~IX.12]{E}. (These results present the image of the trace on $K_0\left(\cprod{C(X)}{S}\right)$ as the preimage of the eigenvalues of $\left(S^1,R_s\right)$ under the canonical map $\R\ni x\mapsto e^{2\pi i x}\in\T^1$, which are precisely the eigenvalues of the standard suspension of $(X,S)$, cf.\  Lemma~\ref{compareEigenvalues}.) Furthermore, by \cite[Corollary~2.6]{PV}, this trace is actually an isomorphism of $K_0\left(\cprod{C(S^1)}{R_s}\right)$ with $\E\left(Y,T\right)$. Hence, since $\E K\left(S^1,R_s,1\right)\subset K_0\left(\cprod{C(S^1)}{R_s}\right)$,  Theorem~\ref{step2} gives us the desired equality $\E K\left(S^1,R_s,1\right)=K_0\left(\cprod{C(S^1)}{R_s}\right)$.\hfill $\Box$
\medskip

In \cite[Corollary~3.7]{Riedel}, N.\ Riedel proved that given two minimal group rotations $(X_1,S_1)$ and $(X_2,S_2)$, their crossed products $\cprod{C(X_1)}{S_1}$ and
$\cprod{C(X_2)}{S_2}$ are isomorphic if and only if $\E(X_1,S_1)=\E(X_2,S_2)$. 
Since $(X_1,S_1)$ and $(X_2,S_2)$  are conjugate if and only if they have the same eigenvales, cf.\ \cite[Theorem~5.18 and Theorem~5.19]{W},  R.\ Gjerde and \O.\ Johansen \cite[Theorem~6]{GJ1999} strengthened Riedel's result by showing that $(X_1,S_1)$ and $(X_2,S_2)$ are conjugate if and only if they are strong orbit equivalent (cf.\ \cite[Definition~1.3]{GPS}). Now,  T.\ Giordano, I.\ Putnam and C.\ Skau \cite[Theorem~2.1]{GPS} proved that two minimal Cantor systems are  strong orbit equivalent if and only if their crossed products are isomorphic. This might lead one to  suspect that the eigenvalues of minimal Cantor systems are invariant of strong orbit equivalence. That this is not the case is probably known by the expert. However, since we could not find such an example in the literature, we give one in the next. 

\begin{example}
There are two (strong) orbit equivalent dynamical systems $(X_1,S_1)$ and $(X_2,S_2)$ for which $\E(X_1,S_1)\not=\E(X_2,S_2)$.
\end{example}
Let $(X_1,S_1)$ be the Denjoy system of Example~\ref{containedStrictly} and let $(X_2,S_2)$ be the Denjoy system with irrational rotation number $s_2$ and with $Q(S_2)=\{ e^{2\pi ns_2}\colon n\in\Z\}\cup\{e^{2\pi i(s_1+ns_2)}\colon\in\Z\}$. A combination of \cite[Theorem~2.2]{GPS} and \cite[Theorem~5.3]{PSS} lead us to conclude that $(X_1,S_1)$ and $(X_2,S_2)$ are (stong) orbit equivalent. By Example~\ref{Denjoy}, the time-$\ts{\frac{1}{s_1}}$ on the standard suspension of $S_1$ is not minimal. Theorem~\ref{step1} tells us then that we can find a (nonzero) rational   number $r$ such that $rs_1$ is an eigenvalue of this standard suspension. But $rs_1$ is not an eigenvalue for the standard suspension of $S_2$; for if it were, then using again Theorem~\ref{step1},  the time-$\textstyle\frac{1}{s_1}$ map on the standard suspension of $S_2$ would not be minimal and so, by  Example~\ref{Denjoy}, $s_1$ would belong to $\Q+s_2\Q$, contradicting the fact that $1$, $s_1$ and $s_2$ are linearly independent over $\Q$.\hfill$\Box$\medskip

I am grateful to T.\ Katsura and J.\ Packer for the discussions we had after I presented the main results of this paper at the Great Planes Operator Theory Symposium 2004. Our last example below was inspired by these discussions. We first prove the following.

\begin{proposition}\label{Packer-Katsura}
Let $(X,S)$ be a minimal dynamical system with at least to distinct $S$-invariant Borel probability measures.  Assume that all normalized traces on $\cprod{C(X)}{S}$ agree on $K_0\left(\cprod{C(X)}{S}\right)$. Then there exists a strictly positive continuous function $f\colon X\to \R$ such that, for all nonzero real number $t$, 
the time-$t$ map on the suspension with base $S$ and ceiling function $f$ is minimal. 
\end{proposition}
\begin{proof}
Let $\mu_0$ and $\mu_1$ be two distinct $S$-invariant Borel probability measures on $X$. For each $t$ in $[0,1]$, put 
\[
    \mu_t = (1-t)\mu_0 + t\mu_1.
\]
By hyphotesis we have that
\begin{equation}\label{sameTrace}
   \tau_{\mu_t} \left(K_0\left(\cprod{C(X)}{S}\right)\right) = \tau\sb{\mu_0}
                    \left(K_0\left(\cprod{C(X)}{S}\right)\right)
 \end{equation}
 for all $t\in [0,1]$. As $\mu_0\not=\mu_1$, there is a continuous function $f\in C(X)$ such that $\tau_{\mu_0}(f)\not=\tau_{\mu_1}(f)$, cf.\ \cite[Theorem~6.2]{W}. Replacing $f$ by a constant plus the imaginary or real part of $f$, if necessary, we may assume that $f$ is a strictly positive real valued function. Let $(Y,T)$ be the suspension with base $S$ and ceiling $f$.  For each $t$ in $[0,1]$ we then get:
\begin{align*}
    \E(Y,T) &\subset \frac{1}{\tau_{\mu_t}(f)}\ \tau_{\mu_t}\left( 
                       K_0\left(\cprod{C(X)}{S}\right)\right)
                       \text{ by Theorem~\ref{step2}}\\
      &= \frac{1}{\tau_{\mu_t}(f)}\ \tau_{\mu_0}\left(K_0\left(\cprod{C(X)}{S}\right)\right)
                       \text{ by (\ref{sameTrace})}.\\      
\end{align*}
Thus
\[
   \E(Y,T)\subset\bigcap_{t\in[0,1]}\frac{1}{\tau_{\mu_t}(f)}\ \tau_{\mu_0}
                     \left(K_0\left(\cprod{C(X)}{S}\right)\right).
\]
As $K_0\left(\cprod{C(X)}{S}\right)$ is countable and $\{\tau_{\mu_t}(f)\}_{t\in[0,1]}$ takes all possible values between $\tau_{\mu_0}(f)$ and $\tau_{\mu_1}(f)$, the right hand side of the above inclusion must be equal to $\{0\}$. Thus $\E(Y,T)=\{0\}$. 
Using now Proposition~\ref{step1}, we conclude that the time-$t$ map on $Y$ is  minimal for all nonzero real number $t$. 
\end{proof}

\begin{example}
There is a dynamical system $(X,S)$ and there is a strictly positive continuous function $f\colon X\to\R$ such that, for all nonzero real number $t$, the time-$t$ map on the suspension with base $S$ and ceiling $f$ is minimal.
\end{example}
Given $\theta\in [0,1]\setminus\Q$, a continuous function  $\xi\colon S^1\to\R$, and 
$n\in\Z\setminus\{0\}$, we define $S=S_{\theta,\xi}\colon S^1\times S^1\to S^1\times S^1$ to be the inverse of the homeomorphism 
\[
    (z_1,z_2)\mapsto\left(e^{2\pi i\theta} z_1, e^{2\pi i \xi(z_1)}z_{1}^{n}z_{2}\right).
\]
The homeomorphism $S$ is called a Furstenberg transformation. It is known that Furstenberg transformations are minimal (see the remark after the Theorem~2.1 in Section~2.3 of \cite{F}).  It is also know that there are Furstenberg transformations which are not uniquely ergodic (see e.g.\ \cite[p.585]{F}.)  Let $S$ be a Furstenberg transformation which is not uniquely ergodic. In \cite[Example~4.9]{Ph1} it is shown that every normalized trace on $\cprod{C(S^1\times S^1)}{S}$  agree on $K_0\left(\cprod{C(S^1\times  S^1)}{S}\right)$. Hence, by Proposition~\ref{Packer-Katsura}, there is a strictly positive function $f\colon S^1\times S^1\to \R$ such that, for all nonzero real number $t$,  the time-$t$ map on the suspension with base $S$ and ceiling $f$ is minimal. \hfill $\Box$


\flushleft{\textit{\mbox{} \\
Department of Mathematics and Statistics,
University of Ottawa, Canada\\
E-mail: bitzaort@uottawa.ca\\
}}
\end{document}